\documentclass[12pt,leqno]{article}
\usepackage{amssymb}
\usepackage[mathscr]{eucal}
\usepackage{amsmath,amssymb,latexsym,theorem,bbm}
\newfont{\msa}{msam10 scaled\magstep1}
\newfont{\ssmsa}{msam9}
\newfont{\smsa}{msam10}
\newfont{\sms}{msbm10}
\newfont{\sseufb}{eufb9}
\newfont{\seufb}{eufb10}
\newfont{\eufb}{eufb10 scaled\magstep1}
\newfont{\eusb}{eusb10 scaled\magstep1}
\newfont{\hcmr}{cmr17 scaled\magstep5}

\newcommand{\DS}{\displaystyle}

\newcommand{\PP}{\mathsf{P}}

\newcommand{\RR}{\mathbb{R}}

\newcommand{\sSS}{\raise-0.5truemm\hbox{\sms S}}

\newcommand{\cB}{{\mathcal B}}
\newcommand{\cE}{{\mathcal E}}

\newcommand{\cF}{{\cal F}}

\newcommand{\cL}{{\cal L}}

\newcommand{\dd}{\mathsf{d}}

\newcommand{\ee}{\mathsf{e}}

\newcommand{\tV}{\widetilde{V}}

\newcommand{\sleq}{\mbox{\ssmsa\hspace*{0.1mm}\symbol{54}\hspace*{0.1mm}}}
\newcommand{\sgeq}{\mbox{\ssmsa\hspace*{0.1mm}\symbol{62}\hspace*{0.1mm}}}
\renewcommand{\leq}{\mbox{\msa\hspace*{0.9mm}\symbol{54}\hspace*{0.9mm}}}
\renewcommand{\geq}{\mbox{\msa\hspace*{0.9mm}\symbol{62}\hspace*{0.9mm}}}

\newcommand{\bone}{\mathbbm{1}}
\newcommand{\vare}{\varepsilon}

\newcommand{\proofend}{\hfill\mbox{$\Box$}}

\numberwithin{equation}{section}

\theoremstyle{change}
\theorembodyfont{\em}
\newtheorem{Lem}{Lemma.}[section]

\newtheorem{Def}[Lem]{Definition.}
\theorembodyfont{\rm}

\begin{document}

\title{\vspace*{-15mm}
       \bfseries\Large Connection between deriving bridges and radial parts
                        from multidimensional Ornstein-Uhlenbeck processes}

\author{By {\sc M\'aty\'as Barczy} (Debrecen)  and {\sc Gyula Pap} (Debrecen)}

\date{}

\maketitle

\renewcommand{\thefootnote}{}
\vspace*{0.2cm}
 \footnote{The first author has been supported by the Hungarian Scientific
  Research Fund under Grant No.\ OTKA--F046061/2004.}

\vspace*{-15mm}

\begin{abstract}
 First we give a construction of bridges derived from a general Markov process
  using only its transition densities.
 We give sufficient conditions for their existence and uniqueness (in law). 
 Then we prove that the law of the radial part of the bridge with endpoints
  zero derived from a special multidimensional Ornstein-Uhlenbeck process
  equals the law of the bridge with endpoints zero derived from the radial
  part of the same Ornstein-Uhlenbeck process. 
 We also construct bridges derived from general multidimensional
  Ornstein-Uhlenbeck processes.
\end{abstract}

\section{Introduction}

In this paper we are dealing with deriving bridges and radial parts from Markov
 processes.
By a bridge from \ $a$ \ to \ $b$ \ over \ $[0,T]$ \ derived from a Markov
 process \ $Z$ \ we mean a process obtained by conditioning \ $Z$ \ to start in
 \ $a$ \ at time \ $0$ \ and arrive at \ $b$ \ at time \ $T$, \ where \ $T>0$.
\ For the construction of such a bridge we use only transition densities.
Important examples are provided by Wiener bridges and Bessel bridges, which
 have been extensively studied and find numerous applications. 
See, for example, Karlin and Taylor \cite[Chapter 15]{KarTay2}, Fitzsimmons,
 Pitman and Yor \cite{FitPitYor}, Baudoin \cite{Bau}, Privault and Zambrini
 \cite{PriZam} and Yor and Zambotti \cite{YorZam}. 
Our construction of bridges is motivated by Karlin and Taylor \cite{KarTay2}
 and Revuz and Yor \cite{RevYor}.
By the radial part of a process with values in \ $\RR^d$ \ we mean its
 euclidean norm. 

We examine whether the operations deriving bridges and radial parts commute
 starting from the same Markov process. 
In case of a multidimensional standard Wiener process and in case of certain
 multidimensional Ornstein-Uhlenbeck processes we show that the answer is yes
 if we consider bridges with endpoints zero. 
We emphasize that Yor and Zambotti in \cite{YorZam} have already proved this
 for a multidimensional standard Wiener process. 
Moreover, they showed that the law of the radial part of the multidimensional
 Wiener bridge with endpoints different from zero is only equivalent and not
 equal to the law of the corresponding multidimensional Bessel bridge.

We proceed as follows. 
In Section 2 we give a construction of a bridge derived from a general Markov
 process using its transition densities. 
We give sufficient conditions for its existence and uniqueness (in law). 
In Sections 3 and 5 we prove that the operations deriving bridges and radial
 parts commute starting from multidimensional standard Wiener processes and
 from certain multidimensional Ornstein-Uhlenbeck processes, respectively. 
In Section 4 we study bridges derived from general multidimensional
 Ornstein-Uhlenbeck processes.

\section{Construction of bridges}

In what follows, let \ $(E,\cE)$ \ be a complete separable metric space endowed
 with the \ $\sigma$-algebra of its Borel subsets, let \ $T>0$, \ let 
 \ $(Z_t)_{0\sleq t\sleq T}$ \ be a time-homogeneous Markov process with state
 space \ $(E,\cE)$ \ admitting transition densities \ $(p_t^Z)_{0<t\sleq T}$
 \ with respect to a fixed \ $\sigma$-finite measure \ $\lambda$ \ on \ $\cE$
 \ (i.e., \ $\PP(Z_t\in A\,|\,Z_s)=\int_Ap_{t-s}^Z(Z_s,y)\,\lambda(\dd y)$
 \ $\PP$-a.s.\ for all \ $A\in\cE$, \ $0\leq s<t\leq T$), \ and let
 \ $a,b\in E$.

If \ $p_t^Z(x,b)>0$ \ for all \ $x\in E$, \ $0<t\leq T$, \ and
 \begin{align}\label{BRIDGE1}
  p_{s,t}(x,y)
  :=\frac{p_{t-s}^Z(x,y)p_{T-t}^Z(y,b)}{p_{T-s}^Z(x,b)},
  \qquad x,y\in E,\quad 0\leq s<t<T,
 \end{align}
 then by a bridge from \ $a$ \ to \ $b$ \ over \ $[0,T]$ \ derived from \ $Z$
 \ we could understand a Markov process \ $(Y_t)_{0\sleq t\sleq T}$ \ with
 initial distribution \ $\PP(Y_0=a)=1$ \ and with transition densities
 \ $(p_{s,t})_{0\sleq s<t<T}$, \ provided that such a process exists (see,
 e.g., Fitzsimmons, Pitman and Yor \cite[Proposition 1]{FitPitYor}, Fitzsimmons
 \cite[Proposition 2.2]{Fit}).
But this definition does not apply for example in case of a \ $d$-dimensional
 Bessel bridge with \ $d>1$ \ and with \ $b=0$, \ since for the transition
 densities \ $(p_t^R)_{t>0}$ \ of the \ $d$-dimensional Bessel process
 \ $(R_t)_{t\sgeq0}$ \ we have \ $p_t^R(x,0)=0$ \ for all \ $x\geq0$, \ $t>0$
 \ (see, e.g., Revuz and Yor \cite[p.~446]{RevYor} or Section 3).

The motivation how to modify \eqref{BRIDGE1} is inspired by Karlin and Taylor
 \cite[p.~267]{KarTay2} and Revuz and Yor \cite[Chapter XI, \S 3]{RevYor}.
For \ $\vare>0$, \ denote by \ $B(b,\vare)$ \ the open ball in \ $E$ \ with
 centre at \ $b$ \ and radius \ $\vare$. 
\ Let \ $(Y_t^\vare)_{0\sleq t\sleq T}$ \ denote \ $(Z_t)_{0\sleq t\sleq T}$
 \ conditioned that \ $Z_T\in B(b,\vare)$.
\ In virtue of Karlin and Taylor \cite[(9.17)]{KarTay2}, for \ $x,y\in E$,
 \ $0\leq s<t<T$, \ the transition densities of \ $Y^\vare$ \ are given by
 \begin{align*}
  p_{s,t}^{Y^\vare}(x,y)
  =p_{t-s}^Z(x,y)
   \frac{\int_{B(b,\vare)}p_{T-t}^Z(y,z)\,\lambda(\dd z)}
        {\int_{B(b,\vare)}p_{T-s}^Z(x,z)\,\lambda(\dd z)},
 \end{align*} 
 provided that \ $\int_{B(b,\vare)}p_{T-s}^Z(x,z)\,\lambda(\dd z)\not=0$.
\ Indeed, by Proposition 7.2 in Kallenberg \cite{Kal},
 \begin{align*}
  \PP(Y_t^\vare\in A\,|\,Y_s^\vare=x)
  &=\PP(Z_t\in A\,|\,Z_s=x,\,Z_T\in B(b,\vare))
   =\frac{\PP(Z_t\in A,\,Z_T\in B(b,\vare)\,|\,Z_s=x)}
         {\PP(Z_T\in B(b,\vare)\,|\,Z_s=x)}\\
  &=\frac{\int_A\int_{B(b,\vare)}
           p_{t-s}^Z(x,y)p_{T-t}^Z(y,z)\,\lambda(\dd y)\lambda(\dd z)}
         {\int_{B(b,\vare)}p_{T-s}^Z(x,z)\,\lambda(\dd z)}
   =\int_A
     p_{s,t}^{Y^\vare}(x,y)\,
     \lambda(\dd y)
 \end{align*} 
 for all \ $A\in\cE$.
\ We can think of the desired bridge as the limit of \ $Y^\vare$ \ as
 \ $\vare\downarrow0$, \ hence our definition is the following.

\begin{Def}\label{DEF:BRIDGE}
For \ $x,y\in E$ \ and \ $0\leq s<t<T$, \ let
 \begin{align}\label{BRIDGE2}
  p_{s,t}(x,y)
  :=p_{t-s}^Z(x,y)
    \lim_{\vare\downarrow0}
     \frac{\int_{B(b,\vare)}p_{T-t}^Z(y,z)\,\lambda(\dd z)}
          {\int_{B(b,\vare)}p_{T-s}^Z(x,z)\,\lambda(\dd z)}
 \end{align}
 if the right hand side exists, \ and \ $p_{s,t}(x,y):=0$ \ otherwise.

By a bridge from \ $a$ \ to \ $b$ \ over \ $[0,T]$ \ derived from \ $Z$ \ we
 mean a Markov process \ $(Y_t)_{0\sleq t\sleq T}$ \ with initial distribution
 \ $\PP(Y_0=a)=1$, \ with \ $\PP(Y_T=b)=1$ \ and with transition densities
 \ $(p_{s,t})_{0\sleq s<t<T}$ \ provided that such a process exists.
\end{Def}

Note that the Markov process \ $(Y_t)_{0\sleq t\sleq T}$ \ (if it exists) is in
 general not time-homogeneous. 
Moreover, additional conditions on \ $(p_t^Z)_{0<t\sleq T}$ \ are needed to
 assure that \ $(Y_t)_{0\sleq t\sleq T}$ \ admits a version having sample paths
 with some regularity properties such as continuity.

\begin{Lem}\label{exists}
Suppose that \ $(p_{s,t})_{0\sleq s<t<T}$ \ defined by
 \eqref{BRIDGE2} satisfy the following properties:
 \renewcommand{\labelenumi}{{\rm(\roman{enumi})}}
 \begin{enumerate}
  \item for all \ $0\leq s<t<T$, \ the function \ $(x,y)\mapsto p_{s,t}(x,y)$
         \ is measurable,  
  \item for all \ $x\in E$ \ and \ $0\leq s<t<T$, \ the function
         \ $y\mapsto p_{s,t}(x,y)$ \ is a probability density, 
  \item for all \ $x,z\in E$ \ and \ $0\leq s<t<u<T$, \ the Kolmogorov-Chapman
         equation
         \ $p_{s,u}(x,z)=\int_Ep_{s,t}(x,y)p_{t,u}(y,z)\,\lambda(\dd y)$
         \ holds.
 \end{enumerate}
Then there exists a unique probability measure \ $\PP_{a,b,T}^Z$ \ on
 \ $\left(E^{[0,T]},\cE^{[0,T]}\right)$ \ such that the coordinate process
 \ $(X_t)_{0\sleq t\sleq T}$ \ on \ $\left(E^{[0,T]},\cE^{[0,T]}\right)$
 \ under \ $\PP_{a,b,T}^Z$ \ is a bridge from \ $a$ \ to \ $b$ \ over \ $[0,T]$
 \ derived from \ $Z$.

Consequently, if \ $(Y_t)_{0\sleq t\sleq T}$ \ is a bridge from \ $a$ \ to
 \ $b$ \ over \ $[0,T]$ \ derived from \ $Z$ \ then its law on
 \ $\left(E^{[0,T]},\cE^{[0,T]}\right)$ \ is \ $\PP_{a,b,T}^Z$.
\end{Lem}

\noindent
\textbf{Proof.} \ 
For \ $x\in E$, \ $A\in\cE$ \ and \ $0\leq s<t<T$, \ let
 \ $\mu_{s,t}(x,A):=\int_Ap_{s,t}(x,y)\,\lambda(\dd y)$,
 \ $\mu_{s,T}(x,A):=\bone_A(b)$, \ where \ $\bone_A$ \ denotes the indicator
 function of the set \ $A$.
\ Then \ $\mu_{s,t}$ \ is a transition probability for all
 \ $0\leq s<t\leq T$, \ and one can check easily that the Kolmogorov-Chapman
 equation \ $\mu_{s,u}(x,A)=\int_E\mu_{s,t}(x,\dd y)\mu_{t,u}(y,A)$ \ holds for
 all \ $x\in E$, \ $A\in\cE$, \ $0\leq s<t<u\leq T$.
\ By Revuz and Yor \cite[Chapter III, Theorem 1.5]{RevYor}, there exists a
 unique probability measure \ $\PP_{a,b,T}^Z$ \ on
 \ $\left(E^{[0,T]},\cE^{[0,T]}\right)$ \ such that the coordinate process
 \ $(X_t)_{0\sleq t\sleq T}$ \ is Markov under \ $\PP_{a,b,T}^Z$ \ with
 transition probabilities \ $(\mu_{s,t})_{0\sleq s<t\sleq T}$ \ and with
 initial distribution \ $\PP_{a,b,T}^Z(X_0=a)=1$. 

Moreover,
 \ $\PP_{a,b,T}^Z(X_T=b)=\PP_{a,b,T}^Z(X_T=b\,|\,X_0=a)=\mu_{0,T}(a,\{b\})=1$.
\proofend

The proof of the next lemma is trivial.

\begin{Lem}\label{LEMMA:CONT}
If \ $f:E\to\RR$ \ is a continuous function then
 \begin{align}\label{IA} 
  \lim\limits_{\vare\downarrow0}
   \frac{1}{\lambda(B(x,\vare))}\int_{B(x,\vare)}f(z)\,\lambda(\dd z)
  =f(x)
 \end{align}
 for all \ $x\in E$.
\ Consequently, if \ $f,g:E\to\RR$ \ are continuous functions such that
 \ $f(z)=g(z)$ \ $\lambda$-a.e. \ $z\in E$ \ then \ $f(z)=g(z)$ \ for all
 \ $z\in E$.
\end{Lem}

\begin{Lem}\label{LEMMA:KC}
If for each \ $0<t\leq T$, \ the probability density \ $p_t^Z$
 \ satisfies the properties
 \renewcommand{\labelenumi}{{\rm(\roman{enumi})}}
 \begin{enumerate}
  \item the function \ $(x,y)\mapsto p_t^Z(x,y)$ \ is continuous, 
  \item for all \ $x_0\in E$, \ there is a \ $\delta>0$ \ such that
         \ $\sup\limits_{x\in B(x_0,\,\delta)}\,\sup\limits_{y\in E}\,
             p_t^Z(x,y)
            <\infty$,
  \item for all \ $y_0\in E$, \ there is a \ $\delta>0$ \ such that
         \ $\sup\limits_{x\in E}\,\sup\limits_{y\in B(y_0,\,\delta)}\,
             p_t^Z(x,y)<\infty$,
  \item for all \ $y\in E$, \ we have
         \ $\int_Ep_t^Z(x,y)\,\lambda(\dd x)<\infty$,
 \end{enumerate}
 then the Kolmogorov-Chapman equation
 \begin{align}\label{KC}
  p_{s+t}^Z(x,z)
  =\int_Ep_s^Z(x,y)p_t^Z(y,z)\,\lambda(\dd y)
 \end{align}
 holds for all \ $x,z\in E$ \ and all \ $s,t>0$ \ with \ $s+t\leq T$.
\ {\rm(Compare with Fitzsimmons, Pitman and Yor \cite[(2.3)]{FitPitYor},
 Fitzsimmons \cite[(1.9)]{Fit}.)}
\end{Lem}

\noindent
\textbf{Proof.} \ 
For \ $x\in E$, \ $A\in\cE$ \ and \ $0<t\leq T$, \ let
 \ $\mu_t^Z(x,A):=\int_Ap_t^Z(x,y)\,\lambda(\dd y)$.
\ Let us fix \ $s,t>0$ \ with \ $s+t\leq T$.
\ Then for all \ $A\in\cE$ \ the Kolmogorov-Chapman equation
 \ $\mu_{s+t}^Z(x,A)=\int_E\mu_s^Z(x,\dd y)\mu_t^Z(y,A)$ \ holds for
 \ $\PP_{Z_s}$-a.e.\ \ $x\in E$, \ where \ $\PP_{Z_s}$ \ denotes the
 distribution of \ $Z_s$ \ (see, e.g., Kallenberg \cite[Corollary 7.3]{Kal}).
Thus for all \ $A\in\cE$
 \begin{align*}  
  \int_Ap_{s+t}^Z(x,z)\,\lambda(\dd z)
  &=\int_Ep_s^Z(x,y)\left(\int_Ap_t^Z(y,z)\lambda(\dd z)\right)\lambda(\dd y)\\
  &=\int_A\left(\int_Ep_s^Z(x,y)p_t^Z(y,z)\lambda(\dd y)\right)\lambda(\dd z)
   \qquad\text{$\PP_{Z_s}$-a.e.\ \ $x\in E$.} 
 \end{align*}
Hence we obtain that for \ $\PP_{Z_s}$-a.e.\ \ $x\in E$, \ equation \eqref{KC}
 holds for \ $\lambda$-a.e.\ $z\in E$.
\ By assumptions (i) and (iii) and the dominated convergence theorem, both
 sides of equation \eqref{KC} are continuous in \ $z\in E$ \ for every fixed
 \ $x\in E$.
\ By Lemma \ref{LEMMA:CONT}, if \ $x\in E$ \ such that \eqref{KC} holds for
 \ $\lambda$-a.e.\ $z\in E$ \ then it holds for all \ $z\in E$.
\ By assumptions (i), (ii) and (iv) and the dominated convergence theorem, both
 sides of equation \eqref{KC} are continuous in \ $x\in E$ \ for every fixed
 \ $z\in E$.
\ The measure \ $\PP_{Z_s}$ \ is clearly \ $\sigma$-finite, hence, again by
 Lemma \ref{LEMMA:CONT}, we conclude that \eqref{KC} holds for all \ $x,z\in E$
 \ and all \ $s,t>0$ \ with \ $s+t\leq T$.
\proofend

\begin{Lem}\label{cont}
Suppose that the densities \ $(p_t^Z)_{0<t\sleq T}$ \ satisfy the following
 properties:
 \renewcommand{\labelenumi}{{\rm(\roman{enumi})}}
 \begin{enumerate}
  \item for all \ $0<t\leq T$, \ the function \ $(x,y)\mapsto p_t^Z(x,y)$ \ is
         continuous,
  \item for all \ $x,z\in E$ \ and all \ $s,t>0$ \ with \ $s+t\leq T$, \ the
         Kolmogorov-Chapman equation \eqref{KC} holds,
  \item for all \ $x\in E$ \ and all \ $0<t\leq T$, \ we have \ $p_t^Z(x,b)>0$.
  \end{enumerate}
Then \eqref{BRIDGE1} holds, and the functions \ $(p_{s,t})_{0\sleq s<t<T}$
 \ satisfy conditions of Lemma \ref{exists}.
\end{Lem}

\noindent
\textbf{Proof.} \ 
Clearly, assumptions (i), (iii) and Lemma \ref{LEMMA:CONT} imply
 \eqref{BRIDGE1}.
Using (i) and (ii), it is easy to check that the functions
 \ $(p_{s,t})_{0\sleq s<t<T}$ \ satisfy conditions of Lemma \ref{exists}.
\proofend

\begin{Lem}\label{Bessel_bridge}
 Let \ $E=[0,\infty)$, \ let \ $\lambda$ \ be the Lebesgue measure on
 \ $[0,\infty)$, \ and let \ $b=0$.
\ Suppose that the densities \ $(p_t^Z)_{0<t\sleq T}$ \ satisfy the following
 properties:
 \renewcommand{\labelenumi}{{\rm(\roman{enumi})}}
 \begin{enumerate}
  \item for all \ $0<t\leq T$, \ the function \ $(x,y)\mapsto p_t^Z(x,y)$ \ is
         continuous,  
  \item for all \ $x,z\in[0,\infty)$ \ and all \ $s,t>0$ \ with \ $s+t\leq T$,
         \ the Kolmogorov-Chapman equation \eqref{KC} holds, 
  \item for all \ $x,y\in[0,\infty)$ \ and all \ $0\leq s<t<T$, \ the limit
         \ $\lim\limits_{\vare\downarrow0}
            \frac{p_{T-t}^Z(y,\vare)}{p_{T-s}^Z(x,\vare)}$
         \ exists,
  \item for all \ $0\leq s<t<T$ \ and all \ $x\in[0,\infty)$, \ there is a
         \ $\delta>0$ \ such that
         \begin{align*}
          \sup_{y\in[0,\infty)}\,\sup_{0<\vare<\delta}\,
           \frac{p_{T-t}^Z(y,\vare)}{p_{T-s}^Z(x,\vare)}
          <\infty.
         \end{align*}
 \end{enumerate}
Then for all \ $x,y\in[0,\infty)$, \ $0\leq s<t<T$, \ we have
 \begin{align}\label{BRIDGE3}
  p_{s,t}(x,y)
  =p_{t-s}^Z(x,y)
   \lim_{\vare\downarrow0}
    \frac{p_{T-t}^Z(y,\vare)}{p_{T-s}^Z(x,\vare)},
 \end{align}
 and the functions \ $(p_{s,t})_{0\sleq s<t<T}$ \ satisfy conditions of Lemma
 \ref{exists}.
\end{Lem}

\noindent
\textbf{Proof.} \ 
Assumptions (i), (iii) and $\cL$'Hospital's rule yield \eqref{BRIDGE3}.
For every \ $0\leq s<t<T$, \ measurability of \ $(x,y)\mapsto p_{s,t}(x,y)$
 \ follows from \eqref{BRIDGE3} and assumptions (i) and (iii).
For every \ $0\leq s<t<T$ \ and \ $x\in[0,\infty)$, \ the function 
 \ $y\mapsto p_{s,t}(x,y)$ \ is a probability density, since by the assumptions
 and the dominated convergence theorem,
 \begin{align*}
  \int_0^\infty p_{s,t}(x,y)\,\dd y
  &=\int_0^\infty\lim_{\vare\downarrow0}
    \frac{p_{t-s}^Z(x,y)p_{T-t}^Z(y,\vare)}{p_{T-s}^Z(x,\vare)}\,\dd y\\
  &=\lim_{\vare\downarrow0}
     \frac{1}{p_{T-s}^Z(x,\vare)}
     \int_0^\infty p_{t-s}^Z(x,y)p_{T-t}^Z(y,\vare)\,\dd y
   =1.
 \end{align*} 
For every \ $0\leq s<t<u<T$ \ and \ $x,z\in[0,\infty)$, \ the
 Kolmogorov-Chapman equation
 \ $p_{s,u}(x,z)=\int_0^\infty p_{s,t}(x,y)p_{t,u}(y,z)\,\dd y$ \ follows from
 \eqref{BRIDGE3} and assumptions (i)--(iii).
\proofend

\section{The case of a standard $d$-dimensional Wiener process}

Let \ $(B_t)_{t\sgeq 0}$ \ be a standard \ $d$-dimensional Wiener process and
 \ $T>0$ \ be fixed. 
Let \ $(X_t)_{0\sleq t\sleq T}$ \ be the bridge with endpoints zero over
 \ $[0,T]$ \ derived from \ $(B_t)_{t\sgeq 0}$ \ (called the \ $d$-dimensional
 Wiener bridge between \ $0$ \ and \ $0$ \ over \ $[0,T]$).
Let \ $R_t=\|B_t\|$, \ $t\geq0$ \ be the radial part of \ $(B_t)_{t\sgeq 0}$
 \ (called the \ $d$-dimensional Bessel process), where \ $\|\cdot\|$
 \ denotes the euclidean norm. 
Let \ $(Y_t)_{0\sleq t\sleq T}$ \ be the bridge with endpoints zero over
 \ $[0,T]$ \ derived from \ $(R_t)_{t\sgeq 0}$ \ (called the \ $d$-dimensional
 Bessel bridge between \ $0$ \ and \ $0$ \ over \ $[0,T]$). 
\ As it is explained by Yor and Zambotti in \cite{YorZam}, a simple invariance
 by rotation argument implies that the laws of \ $(\|X_t\|)_{0\sleq t\sleq T}$
 \ and \ $(Y_t)_{0\sleq t\sleq T}$ \ coincide.
Intuitively, taking bridges with endpoints zero and taking radial parts
 commutate in case of a standard Wiener process, or, in other words, the radial
 part of a Wiener bridge with endpoints zero is the Bessel bridge with
 endpoints zero.
We want to show this result by computing the transition densities of the
 processes \ $(\|X_t\|)_{0\sleq t\sleq T}$ \ and \ $(Y_t)_{0\sleq t\sleq T}$
 \ to demonstrate our method which will also work for certain multidimensional
 Ornstein-Uhlenbeck processes.

It is well known that the transition densities of the process
 \ $(B_t)_{t\sgeq 0}$ \ is
 $$
  p_t^B(x,y)
  =\frac{1}{(2\pi t)^{d/2}}
   \exp\left\{-\frac{\|x-y\|^2}{2t}\right\},
  \qquad t>0,\quad x,y\in\RR^d.
 $$
To demonstrate how to prove Markov property for the radial part of certain
 Markov processes and how to calculate their transition densities, we consider
 the radial part of \ $(B_t)_{t\sgeq 0}$ \ for \ $d\geq2$.
\ We use the ideas due to Karlin and Taylor
 \cite[Chapter 7, Section 6]{KarTay1} and Revuz and Yor
 \cite[Chapter VI, Proposition 3.1]{RevYor}. 
Taking \ $t>0$, \ $0<t_1<\dots<t_n$, \ $b>0$ \ and
 \ $x^{(1)},\dots,x^{(n-1)},x\in\RR^d$, \ we have
 \begin{align*}
  &P(R_{t_n+t}<b\,|\,B_{t_1}=x^{(1)},\dots,B_{t_{n-1}}=x^{(n-1)},\,B_{t_n}=x)
   =P(R_{t_n+t}<b\,|\,B_{t_n}=x)\\
  &=P(R_t<b\,|\,B_0=x)
   =\int_{\|y\|<b}
     \frac{1}{(2\pi t)^{d/2}}\exp\left\{-\frac{\|x-y\|^2}{2t}\right\}\,\dd y
 \end{align*}
 for almost every \ $x\in\RR^d$ \ (with respect to the Lebesgue measure). 
Introducing polar coordinates \ $y=(y_1,\dots,y_d)$ \ by
 \begin{align*}
  y_1&=r\sin\theta_1\cdots\sin\theta_{d-3}\sin\theta_{d-2}\sin\theta_{d-1},\\
  y_2&=r\sin\theta_1\cdots\sin\theta_{d-3}\sin\theta_{d-2}\cos\theta_{d-1},\\
  y_3&=r\sin\theta_1\cdots\sin\theta_{d-3}\cos\theta_{d-2},\\
  \vdots&\\
  y_{d-1}&=r\sin\theta_1\cos\theta_2,\\
  y_d&=r\cos\theta_1,
 \end{align*}
 we obtain
 \begin{align*} 
  P(R_t<b\,|\,B_0=x)
  =\int_0^b
     \frac{r^{d-1}}{(2\pi t)^{d/2}}
     \exp\left\{-\frac{\|x\|^2+r^2}{2t}\right\}G_d(r,x)\,\dd r
 \end{align*}
 for almost every \ $x\in\RR^d$, \ where
 \begin{align*}
  G_d(r,x)
  =\int_{[0,\pi]^{d-2}\times[0,2\pi]}
    (\sin\theta_1)^{d-2}\cdots(\sin\theta_{d-2})
    \exp\left\{\frac{1}{t}\sum_{k=1}^dx_ky_k\right\}\,
    \dd\theta_1\dots\dd\theta_{d-1}
 \end{align*}
 with \ $x=(x_1,\dots,x_d)$.
\ Clearly, the integral
  \ $\int_{\|y\|<b}\exp\left\{-\frac{\|x-y\|^2}{2t}\right\}\,\dd y$
  \ as a function of \ $x$ \ depends only on \ $\|x\|$, \ hence we may put
 \ $x=(0,\dots,0,\|x\|)$, \ and so we obtain
 \begin{align*}
  P(R_t<b\,|\,B_0=x)
  =\int_0^b
    \frac{r^{d-1}}{(2\pi t)^{d/2}}
    \exp\left\{-\frac{\|x\|^2+r^2}{2t}\right\}H_d(r,x)\,\dd r
 \end{align*} 
 for almost every \ $x\in\RR^d$, \ where
 \begin{align*}
  H_d(r,x)
  =2\pi\int_0^\pi(\sin\theta_1)^{d-2}
   \exp\left\{\frac{r\|x\|}{t}\cos\theta_1\right\}\dd\theta_1\,
   \prod_{k=2}^{d-2}\int_0^\pi(\sin\theta_k)^{d-k-1}\,\dd\theta_k.
 \end{align*}
By Gradstein and Ryzhik \cite[8.431]{GraRyz}, for \ $x\not=0$ \ we have  
 \begin{align*}
  \int_0^\pi(\sin\theta_1)^{d-2}
   \exp\left\{\frac{r\|x\|}{t}\cos\theta_1\right\}\,\dd\theta_1
  =\frac{\Gamma\left(\nu+\frac{1}{2}\right)\Gamma\left(\frac{1}{2}\right)}
        {\left(\frac{r\|x\|}{2t}\right)^\nu}I_\nu\left(\frac{r\|x\|}{t}\right),
 \end{align*}
 where \ $\nu=\frac{d}{2}-1$ \ and \ $I_\nu$ \ denotes the modified Bessel
 function of index \ $\nu$ \ defined by
 $$
  I_\nu(z)=\sum_{m=0}^\infty\frac{(z/2)^{2m+\nu}}{m!\Gamma(\nu+m+1)},
  \qquad z>0.
 $$
Moreover, if \ $k$ \ is a positive integer then
 \ $\int_0^\pi(\sin\theta)^k\,\dd\theta=c_k\frac{(k-1)!!}{k!!}$, \ where
 \ $c_k=\pi$ \ if \ $k$ \ is even and \ $c_k=2$ \ if \ $k$ \ is odd.
Consequently,
 \begin{align*}
  P(R_{t_n+t}<b\,|\,B_{t_1}=x^{(1)},\dots,B_{t_{n-1}}=x^{(n-1)},\,B_{t_n}=x)
  =\int_0^b p_t^R(\|x\|,r)\,\dd r
 \end{align*} 
 for almost every \ $x\in\RR^d$, \ where
 \begin{align}\label{BESSEL}
  p_t^R(x,y)
  =\begin{cases}
    \DS\frac{y^{\nu+1}}{tx^\nu}
       \exp\left\{-\frac{x^2+y^2}{2t}\right\}
       I_\nu\left(\frac{xy}{t}\right)
    &\text{if \ $x,y>0$,}\\[6mm]
    \DS\frac{y^{2\nu+1}}{2^\nu t^{\nu+1}\Gamma(\nu+1)}
       \exp\left\{-\frac{y^2}{2t}\right\}
    &\text{if \ $x=0$, \ $y>0$,}
   \end{cases}
 \end{align}
 and \ $p_t^R(x,0):=\lim\limits_{y\downarrow0}p_t^R(x,y)=0$ \ if \ $x\geq0$.
\ Hence
 \begin{align*}
  P(R_{t_n+t}<b\,|\,B_{t_1},\dots,B_{t_n})
  =\int_0^b p_t^R(R_{t_n},r)\,\dd r
  \qquad\text{$\PP$-a.s.}
 \end{align*} 
Clearly, the process \ $(R_t)_{t\sgeq0}$ \ is adapted to the filtration 
 \ $(\cF_t^B)_{t\sgeq0}$, \ where \ $\cF_t^B:=\sigma(B_s,\,0\leq s\leq t)$,
 \ hence we conclude that \ $(R_t)_{t\sgeq0}$ \ is a time-homogeneous Markov
 process with transition densities \ $(p_t^R)_{t>0}$.
\ Note that formula \eqref{BESSEL} is valid also for \ $d=1$ \ with
 \ $p_t^R(x,0):=\lim\limits_{y\downarrow0}p_t^R(x,y)
               =\sqrt{\frac{2}{\pi t}}\exp\left\{-\frac{x^2}{2t}\right\}$
 \ if \ $x\geq0$ \ (see, e.g., Revuz and Yor \cite[p.~446]{RevYor}).

Obviously, for all \ $t>0$ \ and \ $z\in\RR^d$, \ we have 
 \begin{align*}
  \sup_{x,\,y\in\RR^d}\,p_t^B(x,y)=(2\pi t)^{-d/2},\qquad
  \int_{\RR^d}p_t^B(x,z)\,\dd x=1,
 \end{align*}
 hence by Lemmas \ref{exists}, \ref{LEMMA:KC} and \ref{cont} we obtain the
 existence of the Wiener bridge \ $(X_t)_{0\sleq t\sleq T}$ \ and its
 transition densities
 \begin{align*}
  p_{s,t}^X(x,y)
  =\left(\frac{T-s}{2\pi(t-s)(T-t)}\right)^{d/2}
   \exp\left\{-\frac{\|x-y\|^2}{2(t-s)}-\frac{\|y\|^2}{2(T-t)}
              +\frac{\|x\|^2}{2(T-s)}\right\}
 \end{align*} 
 for all \ $x,y\in\RR^d$ \ and all \ $0\leq s<t<T$.

As in case of the Bessel process, one can prove that
 \ $(\|X_t\|)_{0\sleq t\sleq T}$ \ is again a Markov process and obtain its
 transition densities :
 \begin{align*}
  p_{s,t}^{\|X\|}(x,y)
  =\DS\frac{y^{\nu+1}}{(t-s)x^\nu}\left(\frac{T-s}{T-t}\right)^{\nu+1}
      \exp\left\{-\frac{x^2+y^2}{2(t-s)}-\frac{y^2}{2(T-t)}
                 +\frac{x^2}{2(T-s)}\right\}
      I_\nu\left(\frac{xy}{t-s}\right)
 \end{align*}
 for all \ $0\leq s<t<T$ \ and all \ $x,y>0$, \ and
 \begin{align*}
  p_{s,t}^{\|X\|}(0,y)
  =\DS\frac{y^{2\nu+1}}{2^\nu (t-s)^{\nu+1}\Gamma(\nu+1)}
      \left(\frac{T-s}{T-t}\right)^{\nu+1}
      \exp\left\{-\frac{y^2}{2(t-s)}-\frac{y^2}{2(T-t)}\right\}
 \end{align*}
 for all \ $0\leq s<t<T$ \ and all \ $y>0$. 

The aim of the following discussion is to prove that the densities
 \ $(p_t^R)_{t>0}$ \ satisfy conditions of Lemmas \ref{LEMMA:KC} and
 \ref{Bessel_bridge}.
It is known that
 \begin{align*}
  I_\nu(z)=\frac{(z/2)^\nu}{\Gamma(\nu+1)}[1+O(z^2)]
  \quad\text{as \ $z\downarrow0$,}\qquad
  I_\nu(z)=\frac{\ee^z}{\sqrt{2\pi z}}[1+O(z^{-1})]
  \quad\text{as \ $z\to\infty$.}
 \end{align*} 
(For the second statement see Gradstein and Ryzhik \cite[8.451]{GraRyz}.)
Hence
 \begin{align*}
  c_1\left[z^\nu\bone_{(0,1)}(z)+z^{-1/2}\ee^z\bone_{[1,\infty)}(z)\right]
  \leq I_\nu(z)
  \leq c_2\left[z^\nu\bone_{(0,1)}(z)+z^{-1/2}\ee^z\bone_{[1,\infty)}(z)\right]
 \end{align*} 
 with some \ $0<c_1<c_2$ \ for all \ $z>0$.
\ Thus
 \begin{align}\label{Bessel_estimates}
  c_1 f_t(x,y)
  \leq p_t^R(x,y)
  \leq c_2 f_t(x,y)
 \end{align} 
 for all \ $x,y,t>0$, \ where
 \begin{align*}
  f_t(x,y)
  :=t^{-d/2}y^{d-1}\ee^{-(x^2+y^2)/(2t)}\bone_{(0,1)}(xy/t)
    +t^{-1/2}(y/x)^{(d-1)/2}\ee^{-(x-y)^2/(2t)}\bone_{[1,\infty)}(xy/t).
 \end{align*}  
Using \eqref{Bessel_estimates} we obtain
 \ $\sup\limits_{x\sgeq0}\,\sup\limits_{y\sgeq0}\,p_t^R(x,y)<\infty$ \ for all
 \ $t>0$.
\ Indeed, for all \ $t>0$ \ we have
 \begin{align*}
  &\sup_{0<xy<t}\,p_t^R(x,y)
   \leq c_2t^{-d/2}\sup_{y>0}y^{d-1}\ee^{-y^2/(2t)}
   <\infty,\\
  &\sup_{xy\sgeq t,\,y<x}\,p_t^R(x,y)
   \leq c_2t^{-1/2},\\
  &\sup_{xy\sgeq t,\,y\sgeq x}\,p_t^R(x,y)
   =\sup_{\alpha\sgeq1}\sup_{xy\sgeq t,\,y=\alpha x}\,p_t^R(x,y)
   =\sup_{\alpha\sgeq1}\sup_{x\sgeq\sqrt{t/\alpha}}\,p_t^R(x,\alpha x)\\
  &\leq\sup_{\alpha\sgeq1}\sup_{x\sgeq\sqrt{t/\alpha}}\,
        c_2t^{-1/2}\alpha^{(d-1)/2}\ee^{-(\alpha-1)^2x^2/(2t)}
   =c_2t^{-1/2}\sup_{\alpha\sgeq1}\,
        \alpha^{(d-1)/2}\ee^{-(\alpha-1)^2/(2\alpha)}
   <\infty.
 \end{align*}  
Moreover, for all \ $y,t>0$, \ we have
 \begin{align*}
  \int_0^\infty p_t^R(x,y)\,\dd x
  \leq c_2t^{-d/2}y^{d-1}
       \left(\int_0^{t/y}\ee^{-x^2/(2t)}\,\dd x
             +\int_{t/y}^\infty\ee^{-(x-y)^2/(2t)}\,\dd x\right)
  <\infty.
 \end{align*}  
Furthermore, by \eqref{Bessel_estimates}, for all \ $x>0$ \ and all
 \ $0\leq s<t<T$, \ we have
 \begin{align*}
  \sup_{y>0}\,\sup_{0<\vare<(T-s)/x}\,
   \frac{p_{T-t}^R(y,\vare)}{p_{T-s}^R(x,\vare)}
  \leq\frac{c_2}{c_1}\left(\frac{T-s}{T-t}\right)^{d/2}
      \exp\left\{\frac{x^2}{2(T-s)}+\frac{T-s}{2x^2}\right\}.
 \end{align*} 
Using \ $\lim_{z\to0}z^{-\nu}I_\nu(z)=1/(2^\nu\Gamma(\nu+1))$ \ and Lemmas
 \ref{exists}, \ref{LEMMA:KC}, \ref{Bessel_bridge}, one can prove the existence
 of the Bessel bridge \ $(Y_t)_{0\sleq t\sleq T}$ \ and calculate its
 transition densities.
It turns out that the transition densities of the processes
 \ $(\|X_t\|)_{0\sleq t\sleq T}$ \ and \ $(Y_t)_{0\sleq t\sleq T}$ \ coincide.
By Lemma \ref{exists}, their laws on
 \ $\left([0,\infty)^{[0,T]},(\cB([0,\infty)))^{[0,T]}\right)$ \ coincide.

Note that, as a by-product, we proved that the Kolmogorov-Chapman equation
 \ $\int_0^\infty p_s^R(x,y)p_t^R(y,z)\,\dd y=p_{s+t}^R(x,z)$ \ holds for all
 \ $x,z\geq0$ \ and all \ $s,t>0$, \ hence
 \begin{align}\label{II}
  \int_0^\infty y\ee^{-\gamma y^2}I_\nu(\alpha y)I_\nu(\beta y)\,\dd y
  =\frac{1}{2\gamma}\exp\left\{\frac{\alpha^2+\beta^2}{4\gamma}\right\}
   I_\nu\left(\frac{\alpha\beta}{2\gamma}\right)
 \end{align}
 for all \ $\alpha,\beta,\gamma>0$.
\ (Compare with Gradstein and Ryzhik \cite[8.663]{GraRyz}.)
In other words, we obtained a probabilistic proof of \eqref{II}.

\section{Bridges derived from general multidimensional\\
         Ornstein-Uhlenbeck processes}

Let us consider the \ $d$-dimensional stochastic differential equation (SDE)
 \begin{align}\label{GENOUSDE}
  \begin{cases} 
   \dd Z_t=AZ_t\,\dd t+\Sigma\,\dd W_t,\qquad t\geq0,\\
   \phantom{\dd} Z_0=0,
  \end{cases} 
 \end{align}
 where \ $A\in\RR^{d\times d}$, \ $\Sigma\in\RR^{d\times r}$ \ 
and \ $(W_t)_{t\sgeq 0}$ \ is a standard \ $r$-dimensional Wiener process. 
It is known that there exists a strong solution of equation
 \eqref{GENOUSDE}, \ namely
 \begin{align}\label{GENOU_megoldas}
  Z_t=\int_0^te^{(t-s)A}\Sigma\,\dd W_s,\qquad t\geq 0,
 \end{align}
 and pathwise uniqueness for \eqref{GENOUSDE} holds.
(See, e.g., Karatzas and Shreve \cite[5.6]{KarShr}.) 
The process \ $(Z_t)_{t\sgeq 0}$ \ is a time-homogeneous Gauss-Markov process,
 which is called a general \ $d$-dimensional Ornstein-Uhlenbeck (OU) process.
From \eqref{GENOU_megoldas} we obtain
 \begin{align*}
  Z_t=\ee^{(t-s)A}Z_s+\int_s^t\ee^{(t-u)A}\Sigma\,\dd W_u
 \end{align*}
 for all \ $0\leq s<t$, \ thus the conditional distribution of \ $Z_t$ \ with
 respect to \ $Z_s=x$ \ is a normal distribution with mean \ $\ee^{(t-s)A}x$
 \ and variance matrix
 \begin{align*}
  \int_s^t\ee^{(t-u)A}\Sigma\Sigma^\top\ee^{(t-u)A^\top}\,\dd u
  =\int_0^{t-s}\ee^{(t-s-v)A}\Sigma\Sigma^\top\ee^{(t-s-v)A^\top}\,\dd v.
 \end{align*} 
Hence if \ $\Sigma\Sigma^\top$ \ is a (strictly) positive definite matrix
 (necessarily \ $r\geq d$\,) then \ $(Z_t)_{t\sgeq0}$ \ has transition
 densities \ $(p_t^Z)_{t>0}$ \ given by
 \begin{align}\label{GENOUDEN}
  p_t^Z(x,y)
  =\frac{1}{\sqrt{(2\pi)^d\det(V_t)}}
   \exp\left\{-\frac{1}{2}(y-\ee^{tA}x)^\top V_t^{-1}(y-\ee^{tA}x)\right\}
 \end{align}  
 for all \ $x,y\in\RR^d$ \ and all $t>0$, \ where
 \begin{align*}
  V_t:=\int_0^t\ee^{(t-v)A}\Sigma\Sigma^\top\ee^{(t-v)A^\top}\,\dd v,
  \qquad t>0.
 \end{align*}  
We also have
 \begin{align*}
  p_t^Z(x,y)
  =\frac{1}{\sqrt{(2\pi)^d\det(V_t)}}
   \exp\left\{-\frac{1}{2}(x-\ee^{-tA}y)^\top\tV_t^{-1}(x-\ee^{-tA}y)\right\}
 \end{align*}  
 for all \ $x,y\in\RR^d$ \ and all $t>0$, \ where
 \begin{align*}
  \tV_t:=\int_0^t\ee^{-vA}\Sigma\Sigma^\top\ee^{-vA^\top}\,\dd v,
  \qquad t>0.
 \end{align*}   
If all the eigenvalues of \ $A$ \ have negative real parts then
 \ $V_t=V-\ee^{tA}V\ee^{tA^\top}$, \ $t>0$, \ where \ $V$ \ is the unique
 solution of the algebraic matrix equation \ $AV+VA^\top=-\Sigma\Sigma^\top$
 \ given by \ $V=\int_0^\infty\ee^{uA}\Sigma\Sigma^\top\ee^{uA^\top}\,\dd u$.
\ (See, e.g., Karatzas and Shreve \cite[5.6 A]{KarShr}.) 

Obviously, for all \ $t>0$ \ and \ $z\in\RR^d$, \ we have 
 \begin{align*}
  \sup_{x,\,y\in\RR^d}\,p_t^Z(x,y)=\frac{1}{\sqrt{(2\pi)^d\det(V_t)}},\qquad
  \int_{\RR^d}p_t^Z(x,z)\,\dd x=\det(\ee^{-tA}).
 \end{align*}
Hence by Lemmas \ref{exists}, \ref{LEMMA:KC} and \ref{cont} we obtain the
 existence of the general Ornstein-Uhlenbeck bridge \ $(X_t)_{0\sleq t\sleq T}$
 \ over \ $[0,T]$ \ with endpoints zero and its transition densities
 \begin{align}\label{GENOUBRIDGE}
  \begin{split}
   p_{s,t}^X(x,y)
   &=\sqrt{\frac{\det(\tV_{T-s})}{(2\pi)^d\det(\tV_{t-s}\tV_{T-t})}}\\
   &\phantom{\:}
    \times\exp\bigg\{
           -\frac{1}{2}(x-\ee^{-(t-s)A}y)^\top\tV_{t-s}^{-1}(x-\ee^{-(t-s)A}y)
           -\frac{1}{2}y^\top\tV_{T-t}^{-1}y
           +\frac{1}{2}x^\top\tV_{T-s}^{-1}x\bigg\}
  \end{split}
 \end{align} 
 for all \ $x,y\in\RR^d$ \ and all \ $0\leq s<t<T$.

\section {The case of certain Ornstein-Uhlenbeck processes}

Let us consider the \ $d$-dimensional SDE
 \begin{align}\label{OU_egyenlet}
  \begin{cases} 
   \dd Z_t=aZ_t\,\dd t+\sigma\,\dd W_t,\qquad t\geq0,\\
   \phantom{\dd} Z_0=0,
  \end{cases} 
 \end{align}
 where \ $a,\sigma\in\RR$ \ such that \ $\sigma\not=0$, and
 \ $(W_t)_{t\sgeq 0}$ \ is a standard \ $d$-dimensional Wiener process. 
By \eqref{GENOU_megoldas}, the SDE \eqref{OU_egyenlet} has a strong solution
 given by
 \begin{align}\label{OU_megoldas}
  Z_t=\sigma\int_0^te^{a(t-s)}\,\dd W_s,\qquad t\geq 0,
 \end{align}
 and pathwise uniqueness for \eqref{OU_egyenlet} holds.

Let \ $T>0$ \ be fixed.
Let \ $(X_t)_{0\sleq t\sleq T}$ \ be the bridge with endpoints zero over
 \ $[0,T]$ \ derived from \ $(Z_t)_{t\sgeq 0}$.
\ Let \ $R_t:=\|Z_t\|$, \ $t\geq0$ \ be the radial part of
 \ $(Z_t)_{t\sgeq 0}$. 
\ Let \ $(Y_t)_{0\sleq t\sleq T}$ \ be the bridge with endpoints zero over
 \ $[0,T]$ \ derived from \ $(R_t)_{t\sgeq 0}$. 
\ Our aim is to show that the transition densities of the processes
 \ $(\|X_t\|)_{0\sleq t\sleq T}$ \ and \ $(Y_t)_{0\sleq t\sleq T}$ \ coincide.
In fact, we obtain the result of Section 3 as a special case with \ $a=0$ \ and
 \ $\sigma=1$.

From \eqref{GENOUDEN} we obtain the transition densities of the OU process
 \ $(Z_t)_{t\sgeq 0}$ :
 $$
  p_t^Z(x,y)
  =\frac{1}{(2\pi\sigma^2\kappa(a,t))^{d/2}}
   \exp\left\{-\frac{\|y-\ee^{at}x\|^2}{2\sigma^2\kappa(a,t)}\right\},
  \qquad t>0,\quad x,y\in\RR^d,
 $$ 
 where \ $\kappa(a,t)=\frac{\ee^{2at}-1}{2a}$ \ for \ $a\not=0$, \ and
 \ $\kappa(0,t)=t$.
 
As in case of the \ $d$-dimensional Bessel process, one can prove that
 \ $(R_t)_{t\sgeq0}$ \ is a time-homogeneous Markov process with transition
 densities
 \begin{align*}
  p_t^R(x,y)
  =\begin{cases}
    \DS\frac{\ee^{-a\nu t}y^{\nu+1}}{\sigma^2\kappa(a,t)x^\nu}
       \exp\left\{-\frac{\ee^{2at}x^2+y^2}{2\sigma^2\kappa(a,t)}\right\}
       I_\nu\left(\frac{\ee^{at}xy}{\sigma^2\kappa(a,t)}\right)
    &\text{if \ $x,y>0$,}\\[6mm]
    \DS\frac{y^{2\nu+1}}
            {2^\nu(\sigma^2\kappa(a,t))^{\nu+1}\Gamma(\nu+1)}
       \exp\left\{-\frac{y^2}{2\sigma^2\kappa(a,t)}\right\}
    &\text{if \ $x=0$, \ $y>0$,}
   \end{cases}
 \end{align*}
 where \ $\nu=\frac{d}{2}-1$, \ and
 \ $p_t^R(x,0):=\lim\limits_{y\downarrow0}p_t^R(x,y)=0$ \ if \ $d\geq2$,
 \ $x\geq0$, \ and
 \ $p_t^R(x,0):=\lim\limits_{y\downarrow0}p_t^R(x,y)
               =\sqrt{\frac{2}{\pi\sigma^2\kappa(a,t)}}
                \exp\left\{-\frac{\ee^{2at}x^2}{2\sigma^2\kappa(a,t)}\right\}$
 \ if \ $d=1$, \ $x\geq0$.

By \eqref{GENOUBRIDGE}, the transition densities of the OU bridge
 \ $(X_t)_{0\sleq t\sleq T}$ \ is
 \begin{align*}
  p_{s,t}^X(x,y)
  &=\left(\frac{\kappa(a,T-s)}
               {2\pi\sigma^2\kappa(a,t-s)\kappa(a,T-t)}\right)^{d/2}\\[6mm]
  &\phantom{=\:}\DS\times
   \exp\left\{-\frac{\|y-\ee^{a(t-s)}x\|^2}{2\sigma^2\kappa(a,t-s)}
              -\frac{\ee^{2a(T-t)}\|y\|^2}{2\sigma^2\kappa(a,T-t)}
              +\frac{\ee^{2a(T-s)}\|x\|^2}{2\sigma^2\kappa(a,T-s)}\right\}
 \end{align*} 
 for all \ $0\leq s<t<T$ \ and all \ $x,y\in\RR^d$.

As in case of the \ $d$-dimensional Bessel process, one can prove that
 \ $(\|X_t\|)_{0\sleq t\sleq T}$ \ is again a Markov process
 and one can calculate its transition densities :
 \begin{align*}
  p_{s,t}^{\|X\|}(x,y)
  &=\DS\frac{\ee^{-a\nu(t-s)}y^{\nu+1}}{\sigma^2\kappa(a,t-s)x^\nu}
      \left(\frac{\kappa(a,T-s)}{\kappa(a,T-t)}\right)^{\nu+1}
      I_\nu\left(\frac{\ee^{a(t-s)}xy}{\sigma^2\kappa(a,t-s)}\right)\\[6mm]
  &\phantom{=\:}\DS\times
      \exp\left\{-\frac{\ee^{2a(t-s)}x^2+y^2}{2\sigma^2\kappa(a,t-s)}
                 -\frac{\ee^{2a(T-t)}y^2}{2\sigma^2\kappa(a,T-t)}
                 +\frac{\ee^{2a(T-s)}x^2}{2\sigma^2\kappa(a,T-s)}\right\}
 \end{align*}
 for all \ $0\leq s<t<T$ \ and all \ $x,y>0$, \ and
 \begin{align*}
  p_{s,t}^{\|X\|}(0,y)
  &=\DS\frac{y^{2\nu+1}}
            {2^\nu(\sigma^2\kappa(a,t-s))^{\nu+1}\Gamma(\nu+1)}
      \left(\frac{\kappa(a,T-s)}{\kappa(a,T-t)}\right)^{\nu+1}\\[6mm]
  &\phantom{=\:}\DS\times
      \exp\left\{-\frac{y^2}{2\sigma^2\kappa(a,t-s)}
                 -\frac{y^2}{2\sigma^2\kappa(a,T-t)}\right\}
 \end{align*}
 for all \ $0\leq s<t<T$ \ and all \ $y>0$. 

As in Section 3, one can check that the densities \ $(p_t^R)_{t>0}$ \ satisfy
 conditions of Lemmas \ref{LEMMA:KC} and \ref{Bessel_bridge} and one obtains
 the existence of the bridge \ $(Y_t)_{0\sleq t\sleq T}$ \ and its transition
 densities.
It turns out that the transition densities of the processes
 \ $(\|X_t\|)_{0\sleq t\sleq T}$ \ and \ $(Y_t)_{0\sleq t\sleq T}$ \ coincide.
By Lemma \ref{exists}, their laws on
 \ $\left([0,\infty)^{[0,T]},(\cB([0,\infty)))^{[0,T]}\right)$ \ coincide.

\vskip5mm

\parbox{8cm}{M\'aty\'as Barczy\\
             Faculty of Informatics\\
             University of Debrecen\\
             Pf.12\\
             H--4010 Debrecen\\
             Hungary\\ \\
             barczy@inf.unideb.hu}
\hfill
\parbox{5cm}{Gyula Pap\\
             Faculty of Informatics\\
             University of Debrecen\\
             Pf.12\\
             H--4010 Debrecen\\
             Hungary\\ \\
             papgy@inf.unideb.hu}
\end{document}